\documentclass{amsart}
\usepackage{amsmath}\usepackage{amsthm}
\usepackage{tikz}
\usepackage{latexsym}
\usepackage[psamsfonts]{amssymb}

\newtheorem*{cor*}{Corollary}
\newtheorem{theorem}{Theorem}
\newtheorem*{theorem*}{Theorem}

\newtheorem*{lemma*}{Lemma}
\newtheorem*{conj*}{Conjecture}

\theoremstyle{remark}

\newtheorem*{remark*}{Remark}

\newcommand{\comment}[1]{}

\def\sgn{\operatorname{sgn}}  

 \def\Q{\mathbb{Q}}\def\R{\mathbb{R}}\def\Z{\mathbb{Z}}
\def\le{\leqslant} \def\ge{\geqslant}
\def\Re{\operatorname{Re}}

\def\SL{\mathrm{SL}} \def\PSL{\mathrm{PSL}}\def\GL{\mathrm{GL}}\def\PGL{\mathrm{PGL}}

  \def\M{\mathcal{M}}\def\F{\mathcal F}\def\RR{\mathcal{R}}

  \def\T{\mathcal T}\def\cc{\mathcal{C}}
 \def\DD{\Delta} \def\G{\Gamma}
\def\dd{\delta} 
  
\def\a{\alpha}\def\g{\gamma}
\def\+{\,+\,}   \def\={\;=\;}
\def\sm#1#2#3#4{\left(\begin{smallmatrix} #1 & #2  \\ #3 & #4 
\end{smallmatrix}\right)}
\def\be{\begin{equation}}  \def\ee{\end{equation}}

\def\tt{\bigtriangleup}
\def\hh{\mathfrak{H}}
\def\H{\mathcal{H}}
\def\btz{\begin{tikzpicture}}  \def\etz{\end{tikzpicture}}

\begin{document}

\title[Refinement of the Kronecker-Hurwitz class number relation]{
A combinatorial refinement of the \\ Kronecker-Hurwitz class number relation}

\author{Alexandru A. Popa}
\address{Institute of Mathematics ``Simion Stoilow" of the Romanian Academy,
P.O. Box 1-764, RO-014700 Bucharest, Romania}
\email{alexandru.popa@imar.ro}
\thanks{The first author was partly supported by CNCSIS grant 
TE-2014-4-2077. He would like to thank the MPIM in Bonn and the 
IHES in Bures-sur-Yvette for providing support and a stimulating 
research environment while working on this paper.}

\author{Don Zagier}
\address{Max-Planck-Institut f\"ur Mathematik, Vivatsgasse 7, 53111 Bonn, Germany}
\email{don.zagier@mpim-bonn.mpg.de}
\begin{abstract}
We give a refinement of the Kronecker-Hurwitz class number 
relation, based on a tesselation of the Euclidean plane into semi-infinite 
triangles labeled by $\PSL_2(\Z)$ that may be of independent interest. 
\end{abstract}
\maketitle

\section{A refinement of a classical class number relation}

We give a refinement, and a new proof, of the following classical result \cite{G,H,K}.

\begin{theorem}[Kronecker, Gierster, Hurwitz] \label{T1} For any $n\ge 1$ we have
\[
\sum_{t^2\le 4n} H(4n-t^2) \= \sum_{n=ad\atop a,\,d>0} \max(a,d)\;. 
\]
\end{theorem}
\noindent Here $H(D)$ ($D\ge 0, D\equiv 0,3 $ mod 4) is
the Kronecker-Hurwitz class number, which has initial values
\begin{center}
\begin{tabular}{r|rrrrrrrrrrrrr}
    $D$ & 0 & 3 & 4 & 7 & 8 & 11 & 12 & 15 & 16 & 19& 20 & 23 & 24\\ 
    \hline \vspace{2mm}
  $H(D)$& $-\dfrac 1{12}$ & $\dfrac 13$ & $\dfrac 12$ & 1 &1&1& 
  $\dfrac 43$ & 2 & $\dfrac 32$ & 1& 2& 3 & 2  \\
\end{tabular}
\end{center}
and for $D>0$ equals the number of $\PSL_2(\Z)$-equivalence classes 
of  positive definite integral binary quadratic forms of discriminant~$-D$, 
with those classes that contain a multiple of $x^2+y^2$ or of $x^2-xy+y^2$ 
counted with multiplicity $1/2$ or $1/3$, respectively.

Let $\G=\PSL_2(\Z)$. By the $\G$-equivariant bijection 
$\sm abcd \leftrightarrow cx^2 +(d-a) xy -b y^2$
between integral matrices of determinant~$n$ and trace~$t$ 
and quadratic forms of discriminant $t^2-4n$,  
the class number relation can be written as 
 \be \label{0}
 \sum_{\substack{M\in \M_n \\ M \text{ elliptic} }} \chi(z_M) \= 
 \sum_{n=ad\atop a,\,d>0} \max(a,d)\+ 
 \begin{cases}1/6 & \text{ if $n$ is a square,}\\
	\phantom{1}0 & \text{ otherwise,}
 \end{cases}
   \ee
where $\M_n$ is the set of integral matrices of determinant $n$ modulo 
$\pm 1$, $z_M$ is the fixed point of an elliptic $M$ in the upper half-plane $\hh$, 
and $\chi:\hh\rightarrow \Q$ is the modified characteristic function 
of the standard fundamental domain 
 \[ \F \= \{z\in \hh: -1/2\le \Re(z)\le 1/2, \ |z|\ge 1 \} \]
of $\G$ acting on $\hh$ such that $\chi(z)$ is $1/2\pi$ times 
the angle subtended by $\F$ at $z$ (so~$\chi$ is 1 in the interior of $\F$, 
0 outside of $\F$, 1/2 on the boundary points different from the 
corners~$\rho=e^{\pi i/3}$ and~$\rho^2$, and 1/6 at the corners). 
 
We will prove a refinement of~\eqref{0} saying that the subsum of the expression
on the left over all~$M$ in a given orbit of the right action of~$\G$ on~$\M_n$
always takes on one of the values 0, 1, 2 (or~7/6 for the orbit~$\sqrt n\,\G$
if~$n$ is a square).  Specifically, let us define for any right 
coset $K$ in $\M_n/ \G$ (more precisely, $K$ is a right coset in $\PGL_2(\Q)/\G$,
since $\M_n$ is not a group) two positive integers $\delta_K$ and $\delta_K'$
by $\dd_K=\gcd(c,d)$, $\dd_K'=n/\dd_K$, where $\sm abcd$ is any representative of~$K$.  
Then we have:
\begin{theorem}\label{T2} For each right coset $K \in \M_n/ \G$ we have
 \[
 \sum_{\substack{M\in K \\ M \text{ elliptic} }} \chi(z_M)\= 
 1+\sgn(\dd_K'-\dd_K) 
 +\begin{cases} 1/6 & \text{ if $K=\sqrt{n}\; \G$,} \\
     \phantom{1}0 & \text{ otherwise.}  \end{cases}
 \]
\end{theorem}
\noindent Equation~\eqref{0} follows immediately 
by summing the relations in Theorem~\ref{T2} over all cosets in the 
disjoint decomposition $\M_n= \bigsqcup \sm {\dd'}{\beta}0{\dd} \G$ 
with $n=\dd' \dd$ and $0\le \beta<\dd'$.

Theorem~\ref{T2} provides a correspondence between right cosets  
and $\G$-conjugacy classes in~$\M_n$, which generically assigns 
two conjugacy classes to each coset with $\dd'>\dd$. We will deduce it
from a similar statement, Theorem~\ref{T21}, which is sharper in two 
respects (it counts the number of matrices with a fixed point in a
smaller domain, and it allows real coefficients), 
and which gives a generically one-to-one correspondence between cosets 
and conjugacy classes. To state it, we consider a half-fundamental domain  
 \[ \F^- \= \{z\in \hh: -1/2\le \Re(z)\le 0, \ |z|\ge 1 \}\;, \]
and define a function $\a: \GL_2^+(\R) \rightarrow \Q$ by 
\[\a(M) \= \begin{cases} \phantom{1}\chi^{-}(z_M) & 
        \text{ if $M$ is elliptic with fixed point $z_M\in\hh$,  }\\
            \phantom0-\displaystyle\frac 1{12} & \text{ if $M$ is scalar,} \\
              \phantom{\chi M}0 & \text{ if $M$ is parabolic or hyperbolic, }
          \end{cases} \]
where $\chi^{-}$ is defined in the same way as $\chi$ (and hence equals~1 
in the interior of~$\F^-$, 0~outside~$\F^-$, 1/2 on the internal boundary
points of~$\F^-$, and 1/4 and 1/6 at the corners $i$~and~$\rho^2$, respectively). Note that $\a(- M)=\a(M)$, so $\a$ is well-defined on~$M\G$.  
\begin{theorem}\label{T21}
   For $M= \sm yx01\in \GL_2(\R) $ with $y>0$, we have 
 \be \label{1} \sum_{\g\in\G}\a\big(M\g\big)\= 
 \frac{1+\sgn(y-1)}{2} \;. 
 \ee   
\end{theorem}
\noindent Since each coset $K\in \M_n /\G$ 
contains a representative $M$ with $M\infty=\infty$, Theorem~\ref{T2} 
immediately follows from~\eqref{1}, and the fact that the map 
$\pm\sm abcd \mapsto \pm \sm {-a}{b}{c}{-d}$ 
is a bijection between the sets of elements in $\M_n$ having fixed 
point in the left half and in the right half of the standard 
fundamental domain for~$\G$.

Theorem~\ref{T21} is proved in Section~\ref{S3}, as an easy consequence
of a triangulation of a Euclidean half-plane by triangles associated to 
elements of $\G$ (Theorem~\ref{T3}). This triangulation may be of independent
interest, and we give a self-contained treatment in the next section. 

\section{A triangulation of a Euclidean half-plane}   \label{S2}       

Let $\G_\infty=\{\g\in \G \mid g\infty = \infty \}$. We identify 
$\G \smallsetminus \G_\infty$ with a subset of $\SL_2(\Z)$ by 
choosing representatives $\g=\sm abcd$ with $c>0$, and for such~$\g$ 
we define a semi-infinite triangle 
  \be\label{2} \DD(\g)\=\{ (x,y)\in\mathbb R^2 \mid  0\le d-cx-ay\le c\le -dx-by\}\;.
  \ee
(The motivation for this definition is that  $(x,y)\in \DD(\g)$ if and only 
if $\sm yx01 \g$ has a fixed point in $\F^-$.) Note that~$\DD(\g)$ is contained 
in the half-plane
 \[ \H \= \{(x,y)\in\mathbb R^2 \mid y\ge 1\}\; , \]
since $y=c(-dx-by)+d^2-d(d-cx-ay)\ge c^2+d^2-c|d|\ge 1$. 
\begin{theorem} \label{T3} 
We have a tesselation
 \[ \H \;= \bigcup_{\g\in \G\smallsetminus \G_\infty} \DD(\g)  \]
of the half-plane $\H$ into semi-infinite triangles with disjoint interiors.
\end{theorem}
\begin{remark*}
We can extend the triangulation of Theorem \ref{T3} to a triangulation
of all of $\R^2$ by triangles labeled by all of $\G$ if we define
$\DD(\g)$ also for $\g\in\G_\infty$ by 
\[ \DD\big( \sm 1n01 \big) \= [-n-1,-n]\times (-\infty, 1]\;,\]
and can then interpret the extended triangulation as giving a piecewise-linear 
action of~$\G$ on~$\R^2$, with each triangle being a fundamental domain. 
However we will not use this in the sequel. 
\end{remark*}
\begin{proof} The group $\G$ is a free product of its two subgroups generated  
by the elements $S=\sm 0{-1}1{\hfill 0} \text{ and } \ U=\sm 0{-1}1{\hfill 1} $
of orders 2 and 3, respectively,  which fix the two corners of 
$\F^{-}$. Therefore we can view elements of $\G$ as words in 
$S, U, U^2$ or as vertices of the tree shown in Figure~1.
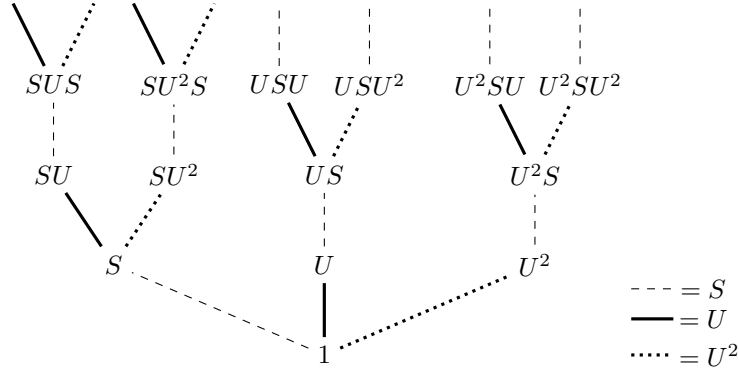
\begin{figure}[h]
\begin{tikzpicture}[scale=0.8,
level 1/.style={sibling distance=35mm},
level 2/.style={sibling distance=20mm},
level 3/.style={sibling distance=15mm},
ss/.style={edge from parent/.style={draw, dashed,thin} },
uu/.style={edge from parent/.style={draw, solid, very thick}},
uu2/.style={edge from parent/.style={draw, dotted, very thick}},
s/.style={dashed, thin},
u/.style={very thick},
u2/.style={very thick, dotted}]

\draw [s] (5.1,1.1)--(5.8,1.1); 
\node at (6.3,1.1) { $= S$}; 
\draw [u] (5.1,0.6)--(5.8,0.6); 
\node at (6.3,0.6) { $= U$}; 
\draw [u2] (5.1,0)--(5.8,0); 
\node at (6.4,0) { $= U^2$}; 

\node at (0,0) {$1$}[grow=up]
 child [uu2]{ node{$U^2$} 
  	child [ss]{ node{$U^2S$}  		
  		child [uu2]{ node{$U^2SU^2$}
  			child[ss]{node{}}  			
  		}
  		child [uu]{ node{$U^2SU$}
  			child[ss]{node{}}  			
  		}
  	}  	
  }
  child[uu]{ node{$U$} 
  	child[ss]{ node{$US$}
  		child[uu2] { node{$USU^2$}
  			child[ss]{node{}}  			
  		}
  		child[uu]{ node{$USU$}
  			child[ss]{node{}}  			
  		}
  	}  	
  }
  child[ss]{ node{$S$}
			child[uu2]{ node{$SU^2$} 
					child[ss]{node{$SU^2S$} 
  							child[uu2]{node{}}  child[uu]{node{}}
  					}
			}
			child[uu]{ node {$SU$}
  					child[ss]{ node{$SUS$} child[uu2]{node{}}
  								child[uu]{node{}}
  					}
  			}
  }; 
\end{tikzpicture}
\caption{A tree associated to $\G= \PSL_2(\Z)$: the vertices 
are labeled by the elements of~$\G$ and the edges by the generators
$S$, $U$ and $U^2$ as shown.}
\label{F0}
\end{figure}
The  proof of both Theorems~\ref{T21}~and~\ref{T3} will follow from 
the following decomposition into triangles with disjoint interiors:
\be \label{3}
   \RR\;:=\:\{(x,y)\in\mathbb R^2  \mid 0\le x\le y-1 \} \= \bigcup_{\g\in \T} \DD(\g) \;, 
\ee
where $\T\subset \G$ is the set of words starting in $U$. 
The regions~$\H$ and $\RR$ and a few triangles corresponding to words of small length are pictured in Figure~\ref{F1}.
\begin{figure}[h]
\begin{tikzpicture}[scale=1.2,
s/.style={dashed, thin},
u/.style={very thick},
u2/.style={very thick, dotted}]
\fill[black!10](0,6)--(0,1)--(5,6)--(0,6);
\draw[->] (0,1) -- (0,6) node[anchor=east] {$y$};
\node at (5.1,0.7) {\footnotesize $y=1$};
\draw (1,2) -- (5,6); 
\draw (1,3) -- (1,6);
\draw (3,5) -- (4,6);
\draw (2,5) -- (2.5,6);
\draw (1,1) -- (4.8,4.8); \node at (1.5,2) {\footnotesize $U^2$};
\draw (2,1) -- (5.5, 2.75); \node at (4.3,2.4) {\footnotesize $U^2SU^2$};
\draw (1,1) -- (5, 3); 
\draw (4,3) -- (5,4); \node at (4.2,3.7) {\footnotesize $U^2SU$};
\draw [u2] (4,3) -- (5,3.5); \draw [s] (3,2)--(4,3);
\node at (5.1,3.2) {\footnotesize $U^2SUS$};
\draw [u] (4,3) -- (5,3.7);

\draw [u](2,1)--(5, 2);\draw (3,1)--(5, 1+2/3);
\draw [s] (5,2)--(5.5,2.25); 

\node at (0.5,2.7) {\footnotesize $U$};
\node at (1.5,2.9) {\footnotesize $U S$};
\node at (1.5,5) {\footnotesize $U SU$};
\node at (2.9,5.7) {\footnotesize $U SU^2$};
\draw (-1,2)--(-1,6); \node at (-0.5,3) {\footnotesize $S$};
\draw (-1,1)--(-4,4); 
\draw (-2,3)--(-3.6,4.6);\node at (-2.3,2.8) {\footnotesize $SU^2$};
\draw (-3,2)--(-4,3);\node at (-2.5 ,2) {\footnotesize $SUS$};
\node at (-1.5 ,3.5) {\footnotesize $SU^2S$};
\draw [s] (-2,5)--(-2.5,6);
\draw [u](-2,3)--(-3,5);\draw (-3,5)--(-3.5,6);
\draw [u2] (-2,3)--(-2,5);\draw (-2,5)--(-2,6);
\node at (-2.8 ,5.6) {\footnotesize $(SU^2)^2$};
\node at (-3.5 ,4.8) {\footnotesize $SU^2SU$};
\draw [s] (-3,5)--(-4,6);
\draw (-2,1)--(-3.6,1.8);\node at (-3.7 ,1.5) {\footnotesize $SUSUS$};
\node at (-3.7 ,2) {\footnotesize $SUSU^2$};
\draw [u2] (-3,1)--(-3.6,1.3);
\draw [u2] (-2,1)--(-3,2);
\draw [s] (-3,2)--(-4,2.5);

\draw [u] (0,1)--(1,2);
\draw [s] (1,2)--(1,3);
\draw [u] (1,3)--(2,5);\draw [s] (2,5)--(2,6);
\draw [u2] (1,3)--(3,5);\draw [s] (3,5)--(3.5,6);

\draw [s] (0,1)--(-4,1);
\draw [u2] (-1,1)--(-1,2);
\draw [s] (-1,2)--(-2,3);

\draw [->,u2] (0,1)--(5.5,1);
\draw [u] (1,1)--(3,2);
\end{tikzpicture}
\caption{The region $\RR$ (shaded) and a few triangles $\tt(\g)$. 
The finite side of a triangle $\DD(\g)$ has been labeled by the final 
letter of~$\g$ as a word in~$S, U, U^2$, with the same convention 
as in Figure~\ref{F0}. }
\label{F1}
\end{figure}
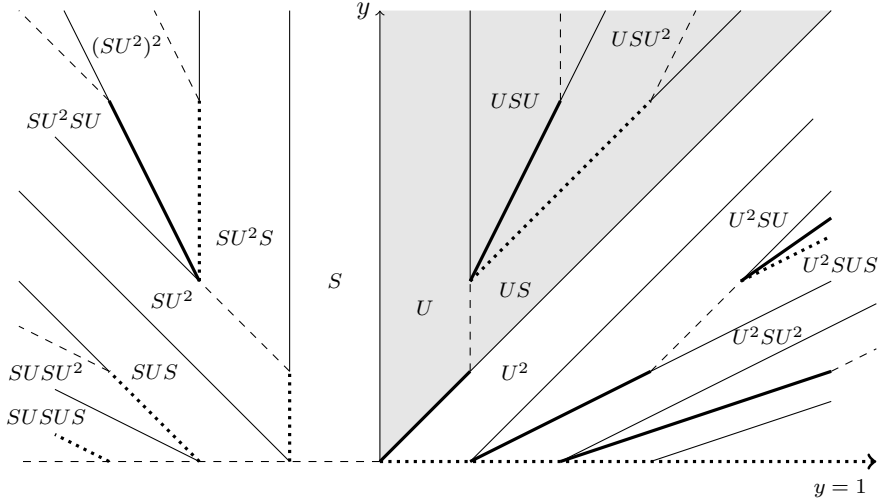

To prove \eqref{3}, let $ \T=\T^+\cup \T^-$, where $\T^+$ consists of the elements of 
$\T$ ending in $U$ or $U^2$, while $\T^-:=\T^+ S$ consists of those elements 
ending in $S$. The set $\T^+$ can be enumerated recursively 
by starting at $U$ and replacing $\g=\sm abcd$ at each step by
\[\g SU= \sm {a}{a+b}{\hfill c}{\hfill c+d}\;,  \quad 
\g SU^2=\sm {a+b}{ b}{\hfill c+d}{\hfill d} \;. \]
From this description we easily obtain the following equivalent 
characterizations\footnote{Recall our convention that $c>0$.}
\[ \g\in\T^+  \iff  0\le \frac {-a}c <\frac {-b}d\le 1, \quad 
   \g\in\T^-  \iff  0\le \frac {-b}d <\frac {-a}c\le 1\;. \]
Alternatively, $\T^+$ consists of those $\g\in\T$ having $d>0$.     

For $\g\in\G\smallsetminus \G_\infty$, the triangle $\tt(\g)$ has two vertices given by
 \[ P_3(-ac-bd+bc,c^2+d^2-cd), \quad  P_2(-ac-bd, c^2+d^2)\;,\]
connected by a line segment of slope $\frac{-d}b$, 
and it has two infinite parallel sides of slope~$\frac{-c}a$. For 
$\g\in\T$ we denote by $\cc(\g)\subset \H$ the 
half-cone containing $\DD(\g)$, bounded by half-lines of slopes 
$-c/a$ and $-b/d$, and having as vertex $P_3$ or $P_2$, depending on
whether $\g\in \T^+$ or $\g\in\T^-$ respectively (see Figure~\ref{F2}). 
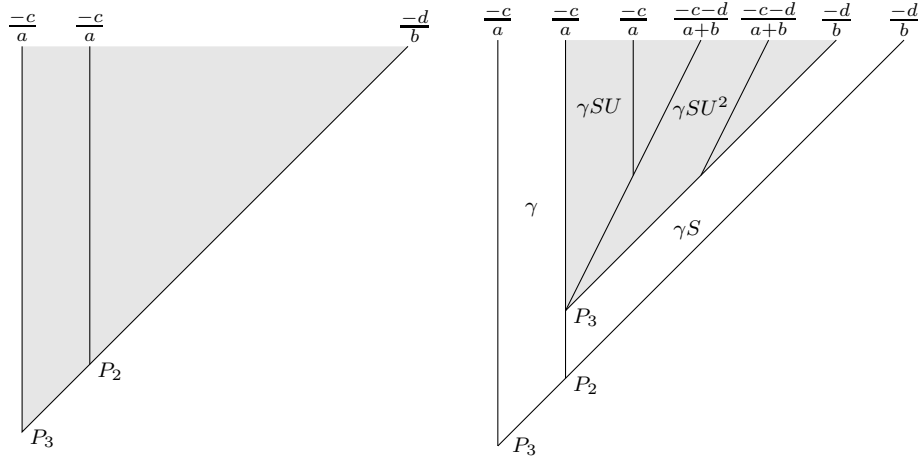
\begin{figure}[h]
\begin{minipage}{.5\textwidth}
\begin{tikzpicture}[scale=0.9]
\fill[black!10](0,6.7)--(0,1)--(5.7,6.7)--(0,6.7);
\draw (0,6.7) -- (0,1)--(5.7,6.7);
\draw (1,2)--(1,6.7);
\node at (0,7) {$\frac{-c}{a}$};
\node at (1,7) {$\frac{-c}{a}$};
\node at (5.8,7) {$\frac{-d}{b}$};
\node at (0.3, 0.9) {\footnotesize $P_3$};
\node at (1.3, 1.9) {\footnotesize $P_2$};
\end{tikzpicture}
\end{minipage}%
\begin{minipage}{.5\textwidth}
\begin{tikzpicture}[scale=0.9]
\fill[black!10](1,7)--(1,3)--(3,7)--(1,7);
\fill[black!10](3,7)--(1,3)--(5,7)--(3,7);
\draw (0,1) -- (0,7);
\draw (0,1) -- (6,7);
\draw (1,2) -- (1,7);
\draw (1,3) -- (5,7);
\draw (1,3) -- (3,7);
\draw (2,5) -- (2,7);
\draw (3,5) -- (4,7);
\node at (0.5,4.5) {\footnotesize $\g$};
\node at (2.8,4.2) {\footnotesize $\g S$};
\node at (1.5,6) {\footnotesize $\g SU$};
\node at (3,6) {\footnotesize $\g SU^2$};
\node at (0,7.3) {$\frac{-c}{a}$};
\node at (1,7.3) {$\frac{-c}{a}$};
\node at (2,7.3) {$\frac{-c}{a}$};
\node at (3,7.3) {$\frac{-c-d}{a+b}$};
\node at (4,7.3) {$\frac{-c-d}{a+b}$};
\node at (5,7.3) {$\frac{-d}{b}$};
\node at (6,7.3) {$\frac{-d}{b}$};
\node at (1.3, 1.9) {\footnotesize $P_2$};
\node at (0.4, 1) {\footnotesize $P_3$};
\node at (1.3, 2.9) {\footnotesize $P_3$};
\end{tikzpicture}
\end{minipage}%
\caption{Left: The cone $\cc(\g)$ and the triangle $\DD(\g)\subset \cc(\g)$  
in the case $\g\in \T^+$. Right: The cone $\cc(\g)$ decomposes into two triangles 
and two smaller, higher-up cones. On top of each line we have marked its 
slope.}\label{F2}
\end{figure}

Using this information, it is easy to check that for 
$\g\in\T^+$ and $\g'=\g S\in \T^-$ we have the 
following decompositions into sets with disjoint interiors (see the right 
picture in Figure~\ref{F2}):
\[
\cc(\g)=\DD(\g) \cup \cc(\g')\;, \quad 
\cc(\g')=\DD(\g')\cup \cc(\g' U)\cup \cc(\g' U^2)\;.
\]
By induction we obtain that $\RR=\cc(U)$ is the union of the triangles
indexed by~$\T$, proving~\eqref{3}.

Finally we show that the decomposition in~\eqref{3} implies the 
decomposition of $\H$ given in Theorem~\ref{T3}. 
From the parenthetical remark following~\eqref{2} it is clear that 
 \[ \DD(T \g)=T \DD(\g)\;, \] 
where $T=SU=\sm 1101$ and $\G_\infty$ acts on $\H$ by  
$T^n(x,y)=(x-ny,y)$. The region 
\be\label{4}
\RR'=\RR\cup \DD(U^2)\=\{(x,y)\in \H: 0\le x< y  \}
\ee
(see Figure~\ref{F1}) is a fundamental domain for this action of $\G_\infty$ 
on $\H$, and we obtain the following decompositions into triangles with disjoint 
interiors
\[\{(x,y)\in\H \mid y-1\le x \}\=\bigcup_{\g\in \T'}\DD(\g),\quad 
  \{(x,y)\in\H \mid x\le 0   \}\=\bigcup_{\g\in \T''}\DD(\g) \;, \]
where $\T'$ consists of words starting with $U^2$, but different 
from $(U^2S)^n=T^{-n}$ with $n>0$, while $\T''$ consists
of words starting with $S$, but different from $(SU)^n=T^n$ with $n>0$. 
Theorem~\ref{T3} follows since 
$\G \smallsetminus \G_\infty = \T \sqcup \T' \sqcup \T''\,$.
\end{proof}

\section{Proof of Theorem \ref{T21}}\label{S3}

Since \eqref{1} is invariant under multiplying $M=\sm yx01$ on the right 
by elements in $\G_\infty$, we assume without loss of generality 
that $0\le x<y \;.$ If $M\g$ is scalar for $\g \in\G$, 
the only possibility is easily seen to be $M=1$. In this case,  
$\a(\g)\ne 0$ for $\g\in \{1, S, U, U^2\}$, 
and~\eqref{1} holds since  $-\frac1{12} + \frac14
+ \frac16 + \frac16 = \frac12$. 

Assuming that $M\ne 1$, it follows that~$\a(M\g)\ne 0$
if and only if~$M\g$ has a fixed point in~$\F^-$, that is 
$(x,y)\in \DD(\g)$. We conclude from Section~\ref{S2} that  $y\ge 1$,
so the point $(x,y)$ belongs to the region  $\RR'$ in~\eqref{4}, and 
$\g=U^2$ or $\g\in\T$ by~\eqref{3}. Therefore the elements $\g$ such 
that~$\a(M\g)\ne 0$ depend on the position of the point $(x,y)$ 
with respect to the triangulation of $\RR'$ as follows (see Figure~\ref{F2}):
\begin{itemize} 
\item $y=1$ and $0<x<1\,$: $\a(MU^2)=1/2$\;;
\item $(x,y)$ is in the interior of a triangle $\tt(\g)\,$: $\a(M\g)=1$\;; 
\item $(x,y)$ is on a common side 
between $\tt(\g)$ and $\tt(\g')$, but it is not a vertex: 
\[\a(M\g)+\a(M\g') \= \frac 12+\frac 12\= 1 \;;\]
\item $(x,y)\in \RR$ is the $P_2$ vertex of the triangle 
$\DD(\g)$ for $\g \in \T^+\,$:
 \[ \a(M\g)+\a(M\g S)+\a(M\g U)\= \frac 14+ \frac 14+\frac 12 \= 1\;;  \]
\item $(x,y)\in \RR$ is the $P_3$ vertex of $\DD(\g')$ with $\g'\in\T^-\,$: 
 \[\a(M\g')+\a(M\g' U)+\a(M\g'U^2)+\a(M\g'S)\= \frac 16+ \frac 16+\frac 
16+\frac 12\= 1\;. 
\]
\end{itemize}

\bibliographystyle{amsplain}

\end{document}